\documentclass[12pt,a4paper]{article}
\usepackage[cp1251]{inputenc}
\usepackage{amsmath}
\usepackage{amssymb}
\usepackage[english]{babel}
\usepackage[dvips]{graphicx}
\usepackage{caption2}
\usepackage{floatflt, curves}
\begin{document}
\title{\textbf{Curvatures of spheres in Hilbert geometry}}
\author{A.A. Borisenko, E.A. Olin}
\maketitle

\begin{center}
\textit{ Geometry Department, Mech.-Math. Faculty, Kharkov National
University, Svoboda sq., 4, 61022-Kharkov, Ukraine.}\\
\textbf{E-mail:} borisenk@univer.kharkov.ua,
evolin@mail.ru
\end{center}

\begin{abstract}
We prove that the normal curvatures of hyperspheres, the Rund curvature, and the Finsler curvature of circles in Hilbert geometry tend to 1 as the radii tend to infinity.
\end{abstract}


\section{Introduction}

A smooth connected manifold $M^n$ is called a \textit{Finsler} manifold (\cite{BCS}) if there is a smooth positively homogeneous on the coordinates in tangent spaces function $F:TM^n\rightarrow [0,\infty)$ such that the symmetric bilinear form
 $\mathbf{g}_y(u,v)=g_{ij}(x,y)u^iv^j:T_x M^n\times T_x M^n\rightarrow
\mathbb{R}$ is positively definite for each pair $(x,y)\in TM^n$, where $g_{ij}(x,y) = \frac{1}{2}[F^2(x,y)]_{y^iy^j}$.

Consider a bounded open convex domain $U$ in $\mathbb{R}^n$ with the Euclidean norm $\|\cdot\|$ and let
$ \partial U$ be a $C^{3}$ hypersurface with positive normal
curvatures. For a point
 $x\in U$ and a tangent vector $y \in
T_xU=\mathbb{R}^n$ let $x_{-}$ and $x_{+}$ be the intersection points of the rays $x+\mathbb{R}_{-}y$ and $x+\mathbb{R}_{+}y$ with \textit{absolute}
$\partial U$. Then the Hilbert metric is defined as follows:
\begin{equation}\label{eq-hilb1}F(x,y)=\frac{1}{2}(\Theta(x,y)+\Theta(x,-y)),\end{equation}
where $\Theta(x,y)=\|y\|\frac{1}{\|x-x_{+}\|}$, $\Theta(x,-y)=\|y\|\frac{1}{\|x-x_{-}\|}$ are called the Funk metrics on $U$.

Hilbert geometries are the generalizations of the Klein's model of the hyperbolic geometry. Hilbert
geometries are also Finsler spaces of constant negative flag curvature $-1$ (\cite{BCS}). Hilbert metric is invariant under projective transformations of $\mathbb{R}^n$ leaving $U$ bounded.

B. Colbois and P. Verovic (\cite{CV}) proved that the Hilbert metric is asymptotically Riemannian at infinity. That means that in a given Hilbert geometry the unit sphere of the norm $F(x,\cdot)$ approaches the ellipsoid in $C^0$ topology as the point $x$ tends to $\partial U$.

Unlike the Riemannian geometry, in the Finsler geometry there are several definitions of the curvature of a curve.

The normal curvature of a hypersurface in a Finsler space is defined as follows (\cite{Sh}).
Let $\varphi:N\rightarrow M^n$ be a hypersurface in a Finsler manifold $M^n$. A vector $\mathbf{n} \in T_{\varphi(x)}M^n$ is called a normal vector to $N$ at the point  $x \in N$ if
$\mathbf{g}_{\mathbf{n}}(y,\mathbf{n})=0$ for all $y \in T_xN$.
The \textit{normal curvature} $\mathbf{k}_\mathbf{n}$ at the point $x\in N$ in a direction $y\in T_xN$ is defined as
\begin{equation}\label{eq-curvform}\mathbf{k}_\mathbf{n} = \mathbf{g}_\mathbf{n}(\nabla_{\dot{c}(s)}\dot{c}(s)|_{s=0},\mathbf{n}),\end{equation}
where $\dot{c}(0)=y$, and $c(s)$ is a geodesic in the induced connection on $N$, $\mathbf{n}$ is the chosen unit normal vector.

For a curve $c(s)$ parameterized by its arc length in $M^n$ it is possible to define two more curvatures.

The Finsler curvature of $c(s)$ (\cite{RR}, \cite{F}) is defined as

\begin{equation}\label{eq-curvf}\mathbf{k}_F(c(s)) = \sqrt{\mathbf{g}_{\dot{c}(s)}(\nabla_{\dot{c}(s)}\dot{c}(s),\nabla_{\dot{c}(s)}\dot{c}(s))} \end{equation}

The Rund curvature of $c(s)$ (\cite{RR}) is defined as
\begin{equation}\label{eq-curvr}\mathbf{k}_R(c(s)) = \sqrt{\mathbf{g}_{\nabla_{\dot{c}(s)}\dot{c}(s)}(\nabla_{\dot{c}(s)}\dot{c}(s),\nabla_{\dot{c}(s)}\dot{c}(s))} \end{equation}

It is well-known that the normal curvatures of hyperspheres in the hyperbolic space $\mathbb{H}^{n}$ are equal to $\coth(r)$ and tend to 1 as the radius  $r$ tends to infinity.
We prove the same property for the Hilbert geometry.

\textbf{Theorem 1.} \textit{The normal curvature, the Rund curvature and the Finsler curvature of the circles centred at the same point in the 2-dimensional Hilbert geometry tend to 1 as their radii tend to infinity, uniformly at the point of the circle.}

\textbf{Theorem 2.} \textit{The normal curvatures of the hyperspheres centred at the same point tend to 1 as their radii tend to infinity, uniformly at the point of the hypersphere and in the tangent vector at this point of the hypersphere.}

That could be interpreted that the Hilbert metric ''tends'' to the Riemannian metric of the hyperbolic space in $C^2$-topology.

\section{The Choice of the Coordinate System}

Consider the Hilbert geometry based on a two-dimensional domain $U$ in the Euclidean plane. Fix a point $o$ in the domain $U$ and a point $p \in \partial U$. Since $\partial U$ is a convex curve, it admits the polar representation $\omega(\varphi)$ from the point $o$ such that the point $p$ corresponds to $\varphi = 0$.

Choose the coordinate system on the plane with the origin $O$ at the point $p$; let the axis $x_2$ be orthogonal to $\partial U$ at $p$, $x_1$ be tangent to $\partial U$ at $p$ and $U-\{p\}$ lie in the half-plane $x_2>0$.

In this section we will construct such a projective transformation $P$ of the plane that sends $U$ to $\hat{U}$ and has the following properties:
\begin{enumerate}
\item $P(p) = p$;
\item The vector $u=(0,1)$ is orthogonal to $\partial \hat{U}$ at the point $p$;
\item The tangent line to $\partial \hat{U}$ at the point $p$ is parallel to the tangent line to $\partial \hat{U}$ at the point corresponding to $\varphi = \pi$;

\item $\partial \hat{U}$ is the graph of the function $x_2=\hat{f}(x_1)$ such that $\hat{f}(0)=0$, $\hat{f}'(0)=0$, $\hat{f}''(0)=\frac{1}{2}$ in the neighbourhood of  $p$;
    \end{enumerate}

 We are going to give the explicit expression for this transformation and show that after this transformation the curvature of $\partial \hat{U}$ and the derivatives of $f$ remain uniformly bounded.

We will use the following lemma that gives the upper bound on the angle between the radial and normal direction to the convex curve.

\textbf{Lemma.} (\cite{Boo2}) \textit{
  Let $\gamma$ be a closed embedded curve in the Euclidean plane whose curvature is greater or equal than $k$. Let $o$ be a point in the interior of the set bounded by $\gamma$, $\omega_0$ be the distance from $o$ to $\gamma$, $\varphi$ be the angle between the outer normal vector at the point $p\in \gamma$ and the vector $op$. Then}

 \begin{equation}
  \label{lem-angle}\cos\angle(u_m,N(m))\geqslant \omega_0 k\end{equation}

Denote by $k$ and $K$ the minimum and maximum of the curvatures of $\partial U$. Denote by $\omega_0=\min\limits_{\varphi}\omega(\varphi)$, $\omega_1=\max\limits_{\varphi}\omega(\varphi)$.

Let the length of the chord of $U$ in the direction $u$ equal $H$, the distance from $o$ to the origin equal $\omega_u$, $\omega_0\leqslant\omega_u\leqslant\omega_1$, and the angle between $u$ and $x_2$ equal $\alpha$.

\textbf{Step 1.} Construct such an affine transformation that makes the vector $\vec{oO}$ parallel to $x_2$. This transformation sends the points $(0,0)$ and $(1,0)$ to themselves, the point $(H \sin \alpha, H \cos \alpha) \in \partial U$ to the point $(0,H)$ and has the expression:
  \begin{equation}
  \label{eq-aff}
           \left\{
           \begin{array}{l}
            \tilde{x}_1= x_1-\tan \alpha x_2 \\
             \tilde{x}_2=\frac{x_2}{\cos \alpha}\\
                    \end{array}
           \right.
       \end{equation}

Denote the image of $U$ as $\tilde{U}$. The point $o$ now has the coordinates $(0,\omega_u)$. Denote by
 $\tilde{k}$ the minimum of the curvature of $\partial \tilde{U}$ in the $(\tilde{x}_1,\tilde{x}_2)$ coordinate system, and by $\tilde{\omega}_0$ denote the distance from the point $(0,\omega_u)$ to $ \partial \tilde{U}$. Note that the eigenvalues of the transformation \eqref{eq-aff} are equal to $1$ and $\frac{1}{\cos \alpha}$, hence
 \begin{equation}
  \label{eq-omega}
    \omega_0 \leqslant \tilde{\omega}_0 \leqslant \frac{1}{\cos \alpha}\omega_0
       \end{equation}

Lemma \eqref{lem-angle} than implies that the curvature of $\partial \tilde{U}$ remains bounded and separated from zero.

\textbf{Step 2.} Construct the transformation such that the tangent line $\tilde{x}_2=-\tan \beta \tilde{x}_1+H$ to $\partial \tilde{U}$ at the point $(0,H)$ will be parallel to the axis~$\tilde{x}_1$, here $\beta$ is the angle between $\tilde{x}_2$ and the normal vector to $\partial \tilde{U}$ at $(0,H)$. This transformation has the expression:
 \begin{equation}
  \label{eq-proj}
           \left\{
           \begin{array}{l}
            \bar{x}_1= \frac{H\tilde{x}_1}{H-\tan \beta\tilde{x}_1} \\
           \bar{x}_2=\frac{H\tilde{x}_2}{H-\tan \beta\tilde{x}_1}\\
                    \end{array}
           \right.
       \end{equation}
Denote the image of $\tilde{U}$ as $\bar{U}$.

  We can estimate the angle $|\tan \beta|$. Using the lemma \eqref{lem-angle} we have,
 \begin{equation}
  \label{eq-tanb}
     0\leqslant |\tan \beta|\leqslant\sqrt{\frac{1}{(\tilde{k}^2 \tilde{\omega}_0^2)}-1}
       \end{equation}

Estimate the curvature $\partial \bar{U}$.

Let the curve $\partial \tilde{U}$  is given in the parametric form $r(t) = (\tilde{x}_1(t),\tilde{x}_2(t))$. Then $\partial \bar{U}$ has the parametrization $\bar{r}(t) = \frac{H r(t)}{H-\tan \beta \tilde{x}_1(t)}$. Differentiating leads to

$$\bar{r}'(t) = \frac{H r'(t)}{H - \tan \beta \tilde{x}_1(t)} + \frac{H r(t) \tan \beta \tilde{x}_1'(t) }{(H - \tan \beta \tilde{x}_1(t))^2} $$

$$ \bar{r}''(t) = \frac{2 H \tan \beta r'(t) \tilde{x}_1'(t)}{(H - \tan \beta \tilde{x}_1(t))^2} + \frac{2 H r(t) \tan^2 \beta \tilde{x}_1'(t)^2 }{(H - \tan \beta \tilde{x}_1(t))^3} +\frac{H r''(t)}{H - \tan \beta \tilde{x}_1(t)} + \frac{H r(t) \tan \beta \tilde{x}_1''(t) }{(H - \tan \beta \tilde{x}_1(t))^2}  $$

The strict convexity of $\partial \tilde{U}$ implies that  $H-\tan \beta \tilde{x}_1(t) \geqslant const > 0$ for each $t$. This and the compactness argument lead to the maximum of the curvature of $\partial \bar{U}$ to be bounded from above for some constant.

If the curve $\partial \tilde{U}$ is the graph $\tilde{x}_2=f(\tilde{x}_1)$, and $f(0)=f'(0) = 0$  then its curvature at the point $(0,0)$ after the transformation \eqref{eq-proj} will not change.
Indeed,

$$\tilde{x}_1'(t)^2+\tilde{x}_2'(t)^2 = $$

$$ \left(\frac{H t \tan \beta}{(H-t \tan \beta)^2}+\frac{H}{H-t \tan \beta}\right)^2+\left(\frac{H \tan \beta f(t)}{(H-t
\tan \beta)^2}+\frac{H f'(t)}{H-t \tan \beta}\right)^2$$

$$\tilde{x}_1'(t)\tilde{x}_2''(t)-\tilde{x}_1''(t)\tilde{x}_2'(t) = $$

$$ = -\left(\frac{2 H t \tan^2 \beta}{(H-t \tan \beta)^3}+\frac{2 H \tan \beta}{(H-t \tan \beta)^2}\right)\left(\frac{H \tan \beta
f(t)}{(H-t \tan \beta)^2}+\frac{H f'(t)}{H-t \tan \beta}\right)+ $$
$$ \left(\frac{H t \tan \beta}{(H-t \tan \beta)^2}+\frac{H}{H-t
\tan \beta}\right)\left(\frac{2 H \tan^2 \beta f(t)}{(H-t \tan \beta)^3}+\frac{2 H \tan \beta f'(t)}{(H-t \tan \beta)^2}+\frac{H
f''(t)}{H-t \tan \beta}\right)$$

We obtain the claim after substituting the equalities $f(0)=f'(0) = 0$.

So, the curvature of $\partial \bar{U}$ at the origin is still separated from zero.

\textbf{Step 3.} Construct such a transformation that the distance from $(0,\omega_u)$ to the origin will equal to 1, and the curvature of $\partial \bar{U}$ at the origin will equal to $1/2$. This transformation has the expression:
 \begin{equation}
  \label{eq-aff2}
           \left\{
           \begin{array}{l}
            \hat{x}_1= \frac{\bar{x}_1}{\omega_u} \\
            \hat{x}_2=\frac{ \bar{x}_2}{2\omega_u^2\bar{k}(0)}\\
                    \end{array}
           \right.
       \end{equation}

Denote the image of $\bar{U}$ as $\hat{U}$. It is obvious that the curvature of  $\partial \hat{U}$ remain bounded.

The announced transformation $P$ is the composition of the transformations \eqref{eq-aff},   \eqref{eq-proj}, \eqref{eq-aff2}, and the following proposition holds:

\textbf{Proposition 1.}\label{th-procurv}\textit{There exists the constant $C_0$ depending on $U$ such that the curvature of $P(\partial U)$ is bounded from above by $C_0$.}

Let $\partial U$ is the graph of the function $x_2=f(x_1)$ in the initial coordinates system. After the transformation P, $P(\partial U)$ can be considered as the graph of the function $x_2=\hat{f}(x_1)$ such that $\hat{f}(0)=0$, $\hat{f}'(0)=0$, $\hat{f}''(0)=\frac{1}{2}$ in the neighborhood of $p$.

Finally estimate the third derivative $\hat{f}'''(0)$. Evidently, under the affine transformations \eqref{eq-aff} and \eqref{eq-aff2} the third derivative remains bounded. We only need to control   $f'''(0)$ at the step 2.

So, let the curve $\partial \tilde{U}$ be the graph $\tilde{x}_2=\tilde{f}(\tilde{x}_1)$ and after the transformation \eqref{eq-proj} we obtain the graph $\bar{f}$. The rules for differentiation lead to
 \begin{equation}\label{eq-3der}\bar{f}'''(0)= \tilde{f}'''(0) - \frac{\tan \beta \tilde{k}(0)}{H} \end{equation}

As $\partial U$ is the compact curve, we obtain

\textbf{Proposition 2.} \textit{There exist the constants $C_1$, $C_2$ depending on $U$, such that $C_1\leqslant\hat{f}'''(0)\leqslant C_2$.}

Analogously we can estimate all higher derivatives.

The Hilbert metrics for the domains $U$ and $\hat{U}$ are isometric. Therefore without loss of generality we will consider the Hilbert metric for the domain $\hat{U}$ and will denote $\hat{U}$  by $U$.

\section{Series Expansions for the Metric Tensor of the Hilbert metric}

From the decomposition of Hilbert metric through the Funk metrics \eqref{eq-hilb1} we conclude
$$
\begin{array}{c}
\displaystyle g_{ij}(x,y)=\mathrm{F}(x,y)\mathrm{F}_{y^iy^j}(x,y)+\mathrm{F}_{y^i}(x,y)\mathrm{F}_{y^j}(x,y)=\\[10pt]
\displaystyle \frac{\mathrm{F}(x,y)(\Theta_{y^iy^j}(x,y)+\Theta_{y^iy^j}(x,-y))}{2}+\\[10pt]
\displaystyle +\frac{(\Theta_{y^i}(x,y)-\Theta_{y^i}(x,-y))(\Theta_{y^j}(x,y)-\Theta_{y^j}(x,-y))}{4}\\[10pt]
\end{array}
$$
Okada lemma (\cite{Sh}) for Funk metrics gives the expression of the derivatives of $\Theta(x,y)$ with respect to the coordinates on tangent spaces through the derivatives with respect to the coordinates on $U$:
$$\Theta(x,y)_{x^k}=\Theta(x,y)\Theta(x,y)_{y^k}$$
Using this lemma we can write:
\begin{equation}
\label{eq-gij}
\begin{array}{c}
\displaystyle g_{ij}(x,y)=\frac{1}{2}\mathrm{F}(x,y)\frac{\Theta_{x^ix^j}(x,y)\Theta(x,y)-2\Theta_{x_i}(x,y)\Theta_{x_j}(x,y)}{\Theta(x,y)^3}+\\[10pt]
\displaystyle +\frac{1}{2}\mathrm{F}(x,y)\frac{\Theta_{x^ix^j}(x,-y)\Theta(x,-y)-2\Theta_{x_i}(x,-y)\Theta_{x_j}(x,-y)}{\Theta(x,-y)^3}+\\[10pt]
\displaystyle +\frac{1}{4}\left(\frac{\Theta_{x^i}(x,y)}{\Theta(x,y)}-\frac{\Theta_{x^i}(x,-y)}{\Theta(x,-y)}\right)\left(\frac{\Theta_{x^j}(x,y)}{\Theta(x,y)}-\frac{\Theta_{x^j}(x,-y)}{\Theta(x,-y)}\right)\\[10pt]
\end{array}
\end{equation}

For convenience we will use lower indices  $x_i$ for coordinates. Let $\mathrm{F}(x_1,x_2,y_1,y_2)$ be a two-dimensional Hilbert metric, $\Theta(x_1,x_2,y_1,y_2)$ be the corresponding Funk metric. Assume that the point $(x_1,x_2)$ is sufficiently close to $\partial U$. Then we can express $\partial U$ as the graph $x_2=f(x_1)$ such that $f(0)=0$, $f'(0)=0$, $f''(0)=\frac{1}{2}$. Consider a point  $(x_1,x_2)$ above the graph $x_2=f(x_1)$. Denote by $r(x_1,x_2,y_1,y_2)$ the distance between the point $(x_1,x_2)$ and the intersection point of the line passing trough $(x_1,x_2)$ in the direction $(y_1,y_2)$ with the curve $x_2=f(x_1)$. Then

 \begin{equation}\label{eq-theta}\Theta(x_1,x_2,y_1,y_2)=\sqrt{y_1^2+y_2^2}r(x_1,x_2,y_1,y_2)^{-1}\end{equation}

Obtain the derivatives of $r(x_1,x_2,y_1,y_2)$ on $x_1$, $x_2$. The parameter $t(x_1,x_2,y_1,y_2)$ corresponding to the intersection points of the curve $x_2=f(x_1)$ with the line

$$
           \left\{
           \begin{array}{rcl}
            x_1(t) = x_1 + t y_1 \\
             x_2(t) = x_2+ t y_2 \\
           \end{array}
           \right.
$$
satisfies the functional equation

\begin{equation}\label{eq-main}x_2+t y_2=f(x_1+ t(x_1,x_2,y_1,y_2) y_1)\end{equation}

Differentiate equation \eqref{eq-main} on $x_1$, $x_2$:

\begin{equation}
\label{eq-mg1}
t_{x_1} y_2=f'(x_1+ t y_1)(1+ t_{x_1}y_1), \ 1+t_{x_2} y_2= f'(x_1+ t y_1)t_{x_2}y_1
\end{equation}

We obtain the explicit expressions for $t_{x_1}$, $t_{x_2}$:

\begin{equation}
\label{eq-mg2}
 t_{x_1}=\frac{f'(x_1+ t y_1)}{y_2-y_1f'(x_1+ t y_1)}, \
 t_{x_2}=\frac{1}{y_1f'(x_1+ t y_1)-y_2}
\end{equation}

Differentiating of \eqref{eq-mg1} leads to:

\begin{equation}
\begin{array}{c}
\label{eq-mg11}
\displaystyle y_2 t_{x_1x_1}=f''(x_1+ t y_1)(1+y_1 t_{x_1})^2+f'(x_1+ t y_1)y_1t_{x_1x_1} \\[10pt]
\displaystyle y_2 t_{x_1x_2}=f''(x_1+ t y_1)(1+y_1t_{x_1})y_2t_{x_2}+f'(x_1+ t y_1)y_1t_{x_1x_2} \\[10pt]
\displaystyle y_2 t_{x_2x_2}=f''(x_1+ t y_1)(y_1 t_{x_2})^2+f'(x_1+ t y_1)y_1t_{x_2x_2} \\[10pt]
 \end{array}
\end{equation}

We obtain the expressions for second derivatives of $t$:

\begin{equation}
\begin{array}{c}
\label{eq-mg22}
\displaystyle t_{x_1x_1}=\frac{f''(x_1+ t y_1)(1+y_1t_{x_1})^2}{y_2-y_1f'(x_1+ t y_1)}\\[10pt]
\displaystyle t_{x_1x_2}=\frac{f''(x_1+ t y_1)(1+y_1t_{x_1})y_1t_{x_2}}{y_2-y_1f'(x_1+ t y_1)}\\[10pt]
\displaystyle t_{x_2x_2}=\frac{f''(x_1+ t y_1)(y_1t_{x_2})^2}{y_2-y_1f'(x_1+ t y_1)}\\[10pt]
\end{array}
\end{equation}

We need the derivatives of $r(x_1,x_2,y_1,y_2)$. By definition:

$$r(x_1,x_2,y_1,y_2)=\sqrt{(y_1 t)^2+(y_2 t)^2}=\sqrt{y_1^2+y_2^2}t(x_1,x_2,y_1,y_2).$$

Hence, $r_{x_k}=t_{x_k}$ and $r_{x_kx_l}=t_{x_kx_l}$.

Now it is possible to calculate the derivatives of the Funk metric. Formula \eqref{eq-theta} implies
\begin{equation}\label{eq-thetader1}\Theta_{x_k}=-\sqrt{y_1^2+y_2^2}\frac{r_{x_k}}{r^2}\end{equation}

After differentiating \eqref{eq-thetader1} we obtain:
\begin{equation}\label{eq-thetader2}\Theta_{x_kx_l}=-\sqrt{y_1^2+y_2^2}\frac{r_{x_kx_l}r^2-2rr_{x_l}r_{x_k}}{r^4}=\sqrt{y_1^2+y_2^2}(2\Theta^3 r_{x_k}r_{x_l}-\Theta^2 r_{x_kx_l})\end{equation}

Finally, from the formula \eqref{eq-gij} it is possible to obtain the coefficients of the metric tensor.

We will need the values of the $g_{ij}(x_1,x_2,y_1,y_2)$ at the points $(x_1,x_2)=(0,x_2)$.

\subsection{Expansions for $g_{ij}(0,x_2,1,0)$}

Note that the strict convexity of $\partial U$ implies that $f'(t(x_1,x_2))\neq 0$ for $t(x_1,x_2) \neq 0 $. Then from \eqref{eq-mg2}  we deduce
\begin{equation}\label{eq-t1}t_{x_1}(0,x_2,1,0)=-1\end{equation}
\begin{equation}\label{eq-t2}t_{x_2}(0,x_2, 1,0)=\frac{1}{f'(t(0,x_2, 1,0))}\end{equation}
and from the formulae \eqref{eq-mg22}:
\begin{equation}\label{eq-t11}t_{x_1x_1}(0,x_2, 1,0)=t_{x_1x_2}(0,x_2, 1,0)=0\end{equation}
\begin{equation}\label{eq-t22}t_{x_2x_2}(0,x_2, 1,0)=-\frac{f''(t(0,x_2, 1,0))}{f'(t(0,x_2, 1,0))^3}\end{equation}
Expanding the functional equation \eqref{eq-main} in a power series with respect to $t$ as $x_2\rightarrow0$ we find the expansions of $t(0,x_2,1,0)$
\begin{equation}\label{eq-tayl}x_2=\frac{1}{4}t^2+\frac{1}{6}f'''(0)t^3+O(t^4)\end{equation}
We will find $t$ in expanded form
\begin{equation}\label{eq-t}t=A+B\sqrt{x_2}+Cx_2+Dx_2^{3/2}+O(x_2^2)\end{equation}
After substituting \eqref{eq-t} into \eqref{eq-tayl} and transposing all members in the left-hand side we obtain the system of equations
\begin{equation}
\label{eq-tt}
\begin{array}{c}
\displaystyle3A^2+2A^3f'''(0)+(6AB+6A^2f'''(0)B)\sqrt{x_2}+\\[10pt]
\displaystyle+(-12+3B^2+6Af'''(0)B^2+6AC+6A^2f'''(0)C)x_2+\\[10pt]
\displaystyle+(2f'''(0)B^2+6BC+12Af'''(0)BC+6AD+6A^2f'''(0)D)x_2^{3/2}+O(x_2^2)=0\\[10pt]
\end{array}
\end{equation}

Choose coefficients $A$, $B$, $C$, $D$ so that the left side of \eqref{eq-tt} is $O(x_2^2)$.
Equating the coefficients under the powers of $x_2$ to zero we obtain two expansions for $t$ which corresponds to the directions $(1,0)$ and $(-1,0)$.
\begin{equation}\label{eq-tpm}t(0,x_2,\pm 1,0)=\pm 2 \sqrt{x_2}-\frac{4}{3}f'''(0)x_2+O(x_2^2)\end{equation}
In our case $r=t$, thus we get
 \begin{equation}\label{eq-rpm}r(0,x_2, 1,0)= 2 \sqrt{x_2}- \frac{4}{3}f'''(0)x_2+O(x_2^2)\end{equation}
 Later on, all power series will be as $x_2\rightarrow0$. Series expansion for the metric $\mathrm{F}$ is:
$$
\begin{array}{c}
\displaystyle \mathrm{F}(0,x_2,1,0)=\frac{1}{2}\left(\frac{1}{r(0,x_2,1,0)}+\frac{1}{r(0,x_2,-1,0)}\right)=\\[10pt]
\displaystyle=\frac{1}{2}\left(\frac{1}{2\sqrt{x_2}- \frac{4}{3}f'''(0)x_2+O(x_2^2)}+\frac{1}{2\sqrt{x_2}+ \frac{4}{3}f'''(0)x_2+O(x_2^2)}\right)=\\[10pt]
\displaystyle=\frac{9}{\sqrt{x_2}(18-8f'''(0)^2x_2)+O(x_2^{3/2})}\\[10pt]
\end{array}
$$
\begin{equation}\label{eq-F}\mathrm{F}(0,x_2,1,0)=\frac{1}{2\sqrt{x_2}}+\frac{2f'''(0)^2}{9}\sqrt{x_2}+O(x_2^{3/2})\end{equation}
We will also need the difference:
\begin{equation}\label{eq-thrazn}
\begin{array}{c}
\displaystyle \Theta(0,x_2,1,0)-\Theta(0,x_2,-1,0)=\frac{1}{r(0,x_2,1,0)}-\frac{1}{r(0,x_2,-1,0)}=\\[10pt]
\displaystyle=\frac{-6f'''(0)+O(x_2)}{4f'''(0)^2x_2-9+O(x_2^{3/2})}=\frac{2}{3}f'''(0)+O(x_2)\\[10pt]
\end{array}
\end{equation}

From \eqref{eq-t1} using $r_{x_k}=t_{x_k}$ we get
\begin{equation}\label{eq-r1}r_{x_1}(0,x_2, 1,0)=-1\end{equation}

Expand the denominator of \eqref{eq-t2} with respect to $t$:
$$
\begin{array}{c}
\displaystyle r_{x_2}(0,x_2, 1,0)=\frac{1}{f'(t(0,x_2, 1,0)}=\\[10pt]
\displaystyle =\frac{1}{f''(0)t(0,x_2, 1,0)+\frac{1}{2}f'''(0)t(0,x_2, 1,0)^2+O(t(0,x_2, 1,0)^3)}\\[10pt]
\end{array}
$$
Using that $f'(0)=0$, $f''(0)=1/2$, and substituting the value of $t$ from \eqref{eq-tpm} we obtain:
$$r_{x_2}(0,x_2,1,0)=\frac{1}{\frac{1}{2}(2 \sqrt{x_2}-\frac{4}{3}f'''(0)x_2)+\frac{1}{2}f'''(0)(2 \sqrt{x_2}-\frac{4}{3}f'''(0)x_2)^2+ O(x_2^2)}$$
And analogously for $r_{x_2}(0,x_2,-1,0)$. Finally,
\begin{equation}\label{eq-r2}r_{x_2}(0,x_2, 1,0)=\frac{1}{\sqrt{x_2}}-\frac{4f'''(0)}{3}+\frac{40f'''(0)^2}{9}\sqrt{x_2}+O(x_2)\end{equation}
The second derivatives has the form $r_{x_kx_l}=t_{x_kx_l}$.
From \eqref{eq-t11} we obtain
\begin{equation}\label{eq-r11}r_{x_1x_1}(0,x_2, 1,0)=r_{x_1x_2}(0,x_2, 1,0)=0\end{equation}
And \eqref{eq-t22} implies
$$r_{x_2x_2}(0,x_2, 1,0)=-\frac{f''(t(0,x_2, 1,0))}{f'(t(0,x_2, 1,0))^3}$$
Expand the numerator and denominator in a series with respect to $t$ and use $f'(0)=0$, $f''(0)=1/2$, \eqref{eq-tpm}:
$$
\begin{array}{c}
\displaystyle r_{x_2x_2}(0,x_2,1,0)=-\frac{\frac{1}{2}+f'''(0)t(0,x_2,1,0)+\frac{1}{2}f^{(4)}(0)t(0,x_2,1,0)^2+O(t^3)}{\left(f''(0)t(0,x_2,1,0)+\frac{1}{2}f'''(0)t(0,x_2,1,0)^2+O(t(0,x_2,1,0)^3)\right)^3}=\\[10pt]
\displaystyle =\frac{-\frac{1}{2}-2f'''(0)\sqrt{x_2}+(\frac{4}{3}f'''(0)^2-4f^{(4)}(0))x_2+\frac{16}{3}f'''(0)f^{(4)}(0)x_2^{3/2}+O(x_2^2)}{x_2^{3/2}+4f'''(0)x_2^2+O(x_2^{5/2})}\\[10pt]
\end{array}
$$

Thus
\begin{equation}\label{eq-r22}r_{x_2x_2}(0,x_2, 1,0)=-\frac{1}{2x_2^{3/2}}-\frac{2f^{(4)}(0)}{ \sqrt{x_2}}+O(1)\end{equation}
From \eqref{eq-thetader1}, \eqref{eq-rpm}, \eqref{eq-r1} we find that
$$\Theta_{x_1}(0,x_2,1,0)=\frac{1}{(2 \sqrt{x_2}- \frac{4}{3}f'''(0)x_2+O(x_2^2))^2}$$

Analogously acting for the vector $(-1,0)$ we get
\begin{equation}\label{eq-theta1}\Theta_{x_1}(0,x_2,\pm1,0)=\pm\frac{1}{4x_2}+\frac{f'''(0)}{3\sqrt{x_2}}+O(1)\end{equation}

From \eqref{eq-rpm}, \eqref{eq-r2} we deduce
$$\Theta_{x_2}(0,x_2,1,0)= -\frac{\frac{1}{\sqrt{x_2}}-\frac{4f'''(0)}{3}+\frac{40f'''(0)}{9}\sqrt{x_2}+O(x_2)}{(2 \sqrt{x_2}- \frac{4}{3}f'''(0)x_2+O(x_2^2))^2}$$

And, finally,
\begin{equation}\label{eq-theta2}\Theta_{x_2}(0,x_2,\pm1,0)=-\frac{1}{4x_2^{3/2}}-\frac{f'''(0)^2}{\sqrt{x_2}}+O(1)\end{equation}

Using the formulae \eqref{eq-thetader2}, \eqref{eq-rpm}, \eqref{eq-r2}, \eqref{eq-r22} we obtain  the expression for the second derivatives of the Funk metric:
\begin{equation}\label{eq-theta22}\Theta_{x_2x_2}(0,x_2,\pm1,0)=\frac{3}{8x_2^{5/2}}+\frac{13f'''(0)^2+3f^{(4)}(0)}{6x_2^{3/2}}+O\left(\frac{1}{x_2}\right)\end{equation}

Finally we can estimate the metric coefficients.  From \eqref{eq-theta}, \eqref{eq-rpm}, \eqref{eq-theta1} we get
\begin{equation}\label{eq-thetabl1}\frac{\Theta_{x_1}(0,x_2,1,0)}{\Theta(0,x_2,1,0)}-\frac{\Theta_{x_1}(0,x_2,-1,0)}{\Theta(0,x_2,-1,0)}=\frac{1}{\sqrt{x_2}}+\frac{4f'''(0)^2}{9}\sqrt{x_2}+O(x_2^{3/2})\end{equation}
It follows from \eqref{eq-theta}, \eqref{eq-rpm}, \eqref{eq-theta2} that
\begin{equation}\label{eq-thetabl2}\frac{\Theta_{x_2}(0,x_2,1,0)}{\Theta(0,x_2,1,0)}-\frac{\Theta_{x_2}(0,x_2,-1,0)}{\Theta(0,x_2,-1,0)}=\frac{2 f'''(0)}{3\sqrt{x_2}}+O(1)\end{equation}
Note that
\begin{equation}\label{eq-thetabl11}
\begin{array}{c}
\displaystyle \Theta_{x_1x_1}(0,x_2,\pm1,0)\Theta(0,x_2,\pm1,0)-2\Theta_{x_1}(0,x_2,\pm1,0)\Theta_{x_1}(0,x_2,\pm1,0)=\\[10pt]
\displaystyle =(2\Theta^3 r_{x_1}r_{x_1}-\Theta^2 r_{x_1x_1})\Theta-2\Theta^2r_{x_1}\Theta^2r_{x_1}=0,\\[10pt]
\end{array}
\end{equation}
since $r_{x_1x_1}=0$. And analogously
\begin{equation}\label{eq-thetabl12}\Theta_{x_1x_2}(0,x_2,\pm1,0)\Theta(0,x_2,\pm1,0)-2\Theta_{x_1}(0,x_2,\pm1,0)\Theta_{x_2}(0,x_2,\pm1,0)=0\end{equation}
Then from \eqref{eq-theta}, \eqref{eq-rpm}, \eqref{eq-theta2}, \eqref{eq-theta22} we get
\begin{equation}
\label{eq-thetabl22}
\begin{array}{c}
\displaystyle\frac{\Theta_{x_2x_2}(0,x_2,\pm1,0)\Theta(0,x_2,\pm1,0)-2\Theta_{x_2}(0,x_2,\pm1,0)\Theta_{x_2}(0,x_2,\pm1,0)}{\Theta(0,x_2,\pm1,0)^3}=\\[10pt]
\displaystyle=\frac{1}{2x_2^{3/2}}+\frac{2f^{(4)}(0)}{\sqrt{x_2}}+O(1)\\[10pt]
\end{array}
\end{equation}
Finally, using \eqref{eq-F}, \eqref{eq-gij}, \eqref{eq-thetabl1}, \eqref{eq-thetabl2}, \eqref{eq-thetabl11}, \eqref{eq-thetabl12}, \eqref{eq-thetabl22} we obtain the series expansions of the metric tensor of Hilbert metric.

\begin{equation}
\label{eq-g11}
\begin{array}{c}
\displaystyle g_{11}(0,x_2,1,0)=\frac{1}{4x_2}+O(1), \  g_{12}(0,x_2,1,0)=\frac{f'''(0)}{6x_2}+O(1)\\[10pt]
\displaystyle g_{22}(0,x_2,1,0)=\frac{1}{4x_2^2}+\frac{2f'''(0)^2+9f^{(4)}(0)}{18x_2}+O(1)\\[10pt]
\end{array}
\end{equation}

\subsection{Expansions for $g_{ij}(0,x_2,0,1)$}

Formulae \eqref{eq-mg2} imply that at the point $(0,x_2)$

$$t_{x_1}(0,x_2,0,\pm 1)=0,\,\, t_{x_2}(0,x_2,0,\pm 1)=-1$$
$$t_{x_1x_2}(0,x_2,0,\pm 1)=t_{x_2x_2}(0,x_2,0,\pm 1)=0$$
Note that the functions $t(0,x_2,0,\pm 1)$ have the representations
$$ t(0,x_2,0,-1)=-x_2,\,\, t(0,x_2,0,1)=H-x_2$$
Here $H$ denotes the length of the chord of $\partial U$ in the direction $(0,1)$.
Then
$$
\Theta(0,x_2,0,-1)=\frac{1}{x_2}, \ \Theta(0,x_2,0,1)=\frac{1}{H-x_2}\\[10pt]
$$
Consequently
\begin{equation}\label{eq-F2}\mathrm{F}(0,x_2,0,1)=\frac{1}{2}\left(\frac{1}{H-x_2}+\frac{1}{x_2}\right)=\frac{1}{2 x_2}+O\left(1\right)\end{equation}
We can estimate the derivatives of the Funk metrics $\Theta(0,x_2,0,\pm1)$. It follows from \eqref{eq-thetader1}, \eqref{eq-thetader2} that
\begin{equation}
\label{eq-thnorm1}
\Theta_{x_2}(0,x_2,0,-1)=\frac{1}{x_2^2}, \ \Theta_{x_2}(0,x_2,0,1)=-\frac{1}{(H-x_2)^2}\\[10pt]
\end{equation}
\begin{equation}
\label{eq-thnorm2}
 \Theta_{x_2x_2}(0,x_2,0,-1)=\frac{2}{x_2^3}, \ \Theta_{x_2x_2}(0,x_2,0,1)=\frac{2}{(H-x_2)^3}\\[10pt]
\end{equation}

Using \eqref{eq-gij}, \eqref{eq-thnorm1}, \eqref{eq-thnorm2}, we get the expansions:
$$g_{12}(0,x_2,0,1)=0$$
\begin{equation}\label{eq-garb}g_{22}(0,x_2,0,1)=\frac{1}{4}\left(\frac{1}{H-x_2}+\frac{1}{x_2}\right)^2=\frac{1}{4 x_2^2}+O\left(\frac{1}{x_2}\right)\end{equation}

We will also need the values $\mathrm{F}(0,x_2,l,\frac{1}{2})$.

We have
$$
\begin{array}{c}
\displaystyle t\left(0,x_2,-l,-\frac{1}{2}\right)=-2x_2+2l^2x_2^2+O(x_2^3)\\[10pt]
\displaystyle t\left(0,x_2,l,\frac{1}{2}\right)=L+O(x_2)\\[10pt]
\end{array}
$$

Then,
$$
\begin{array}{c}
\displaystyle \mathrm{F}\left(0,x_2,l,\frac{1}{2}\right)=\\[10pt]
\displaystyle =\frac{\sqrt{\frac{1}{4}+l^2}}{2\sqrt{\frac{1}{4}t(0,x_2,l,\frac{1}{2})^2+(lt(0,x_2,l,\frac{1}{2}))^2}}+\\[10pt]
\displaystyle +\frac{\sqrt{\frac{1}{4}+l^2}}{2\sqrt{\frac{1}{4}t(0,x_2,-l,-\frac{1}{2})^2+(lt(0,x_2,-l,-\frac{1}{2}))^2}}=\\[10pt]
\displaystyle =\frac{\sqrt{\frac{1}{4}+l^2}}{2\sqrt{\frac{1}{4}+l^2}}\left(\frac{1}{t(0,x_2,l,\frac{1}{2})}-\frac{1}{t(0,x_2,-l,-\frac{1}{2})}\right)\\[10pt]
\end{array}
$$
Finally,
\begin{equation}\label{eq-Fl}\mathrm{F}\left(0,x_2,l,\frac{1}{2}\right)=\frac{1}{4x_2}+\frac{1}{2L}+O(x_2)\end{equation}

\section{Proof of the Theorems}

The Chern-Rund covariant derivative along the curve $c(t)$ in the Finsler space equipped with the Hilbert metric $F$ is given by the formula (\cite{Sh})
\begin{equation}\label{eq-nablacc}\nabla_{c'(t)}c'(t) = \left\{  c''(t)^i+(\Theta(c(t),c'(t))-\Theta(c(t),-c'(t))c'(t)^i \right\} \frac{\partial}{\partial x^i}\end{equation}

For calculating the normal curvature \eqref{eq-curvform}, the Finsler curvature \eqref{eq-curvf} and the Rund curvature  \eqref{eq-curvr} we need the covariant derivative $\nabla_{\dot{c}(s)}\dot{c}(s)$ of the curve $c(s)$ parameterized by its arc length.

For a given curve $c(t)$ we will denote by the dot the derivative with respect to the arc length $s$, and by the prime the derivative with respect to $t$. Then let $t=t(s)$ be the reparameterization. We get
$$\dot{c}(s)=c'(t)t'_s$$
Using that $s$ in the length parameter we get
$$1=\mathrm{F}(c(t),c'(t))t'_s$$
Hence,
$$\dot{c}(s)=\frac{c'(t)}{\mathrm{F}(c(t),c'(t))}$$
Next step is to calculate $\nabla_{\dot{c}(s)}\dot{c}(s)$.

$$
\begin{array}{c}
\displaystyle\nabla_{\dot{c}(s)}\dot{c}(s) = \nabla_{\frac{c'(t)}{\mathrm{F}(c(t),c'(t))}}\frac{c'(t)}{\mathrm{F}(c(t),c'(t))}=\\[10pt]
\displaystyle\frac{1}{\mathrm{F}(c(t),c'(t))}\left( \nabla_{c'(t)}\left(\frac{1}{\mathrm{F}(c(t),c'(t))} \right)c'(t) + \frac{1}{\mathrm{F}(c(t),c'(t))}\nabla_{c'(t)}c'(t)   \right)\\[10pt]
\end{array}
$$
According to (\cite{BCS}),
$$\nabla_{c'(t)}\left(\frac{1}{\mathrm{F}(c(t),c'(t))} \right)=-\frac{\mathbf{g}_{c'(t)}(\nabla_{c'(t)}c'(t),c'(t))}{\mathrm{F}(c(t),c'(t))^3}$$
Then the derivative $\nabla_{\dot{c}(s)}\dot{c}(s)$ has the form
$$\nabla_{\dot{c}(s)}\dot{c}(s)=\frac{1}{\mathrm{F}(c(t),c'(t))^2}\left(\nabla_{c'(t)}c'(t)-\frac{\mathbf{g}_{c'(t)}(\nabla_{c'(t)}c'(t),c'(t))}{\mathrm{F}(c(t),c'(t))^2}c'(t)\right)$$
Finally, using \eqref{eq-nablacc} we get the formula:
\begin{equation}\label{eq-nabl}\nabla_{\dot{c}(s)}\dot{c}(s)=\frac{c''(t)+c'(t)\left(\Theta(c(t),c'(t))-\Theta(c(t),-c'(t))-\frac{\mathbf{g}_{c'(t)}(\nabla_{c'(t)}c'(t),c'(t))}{\mathrm{F}(c(t),c'(t))^2}\right)}{\mathrm{F}(c(t),c'(t))^2}\end{equation}

As in Section 2 fix a point $o$ in the domain $U$ and a point $p \in \partial U$. The curve $\partial U$ admits the polar representation $\omega(\varphi)$ from the point $o$ such that the point $p$ corresponds to $\varphi = 0$. According to Section 2, we assume that $U$ satisfies the conditions 1)-4).

Then one can get that $\omega'(0)=0$, $\omega(0)=1$, $\omega''(0)=1/2$, $\omega'(\pi)=0$. Denote by $C=\frac{1+\omega(\pi)}{\omega(\pi)}$.

In the paper \cite{Boo3} the polar function $\rho_{r}(u)$ of the hypersphere of radius  $r$ was obtained.
\begin{equation}\label{eq-sph}\rho_{r}(u)=\frac{\omega(-u)\omega(u)(e^{2r}-1)}{\omega(u)+\omega(-u)e^{2r}}\end{equation}
as $r\rightarrow\infty$:
\begin{equation}\label{eq-sphas}\omega(u)-\rho_r(u)=\omega(u)\left(\frac{\omega(u)}{\omega(-u)} +
1 \right)e^{-2r}+o(e^{-2r})\end{equation}

From \eqref{eq-sph} we get that the circle of radius $r$ admits the parametrization
$$c(\varphi)=\left(\frac{\omega(\pi-\varphi)\omega(\varphi)(e^{2r}-1)}{\omega(\varphi)+\omega(\pi-\varphi)e^{2r}}\sin \varphi,\frac{\omega(\pi-\varphi)\omega(\varphi)(e^{2r}-1)}{\omega(\varphi)+\omega(\pi-\varphi)e^{2r}}\cos \varphi\right),$$
where $\omega(\varphi)$ is the polar function of $\partial U$.

Then
\begin{equation}\label{eq-c1}c'(0)=\frac{\omega(\pi)(e^{2r}-1)}{1+\omega(\pi)e^{2r}}(1,0)=(1-Ce^{-2r}+O(e^{-3r}),0),r\rightarrow\infty\end{equation}
The second derivative:
$$c''(0)=\frac{(e^{2r}\omega(\pi)^2(\omega''(0)-1)-\omega(\pi)+\omega''(\pi))(e^{2r}-1)}{(1+e^{2r}\omega(\pi))^2}(0,1)$$
\begin{equation}\label{eq-c11}c''(0)=\left(0,-\frac{1}{2}+O(e^{-2r})\right),r\rightarrow\infty\end{equation}
From \eqref{eq-sphas} we get that at the point of the circle the second coordinate
\begin{equation}\label{eq-x2}x_2=\omega(0)-\frac{\omega(\pi)\omega(0)(e^{2r}-1)}{\omega(0)+\omega(\pi)e^{2r}}=C e^{-2r}+O(e^{-3r})\end{equation}

Estimate the derivative $\nabla_{\dot{c}(0)}\dot{c}(0)$ using the formulae \eqref{eq-nabl}, \eqref{eq-thrazn}, \eqref{eq-x2}:
$$
\begin{array}{c}
\displaystyle \Theta(c(0),c'(0))-\Theta(c(0),-c'(0))=\Theta(0,C e^{-2r}+O(e^{-3r}),1+O(e^{-2r}),0)-\\[10pt]
\displaystyle -\Theta(0,C e^{-2r}+O(e^{-3r}),-1+O(e^{-2r}),0)=\frac{2}{3}f'''(0)+ O(e^{-2r})\\[10pt]
\end{array}
$$
Therefore, formula \eqref{eq-nablacc} leads to

\begin{equation}\label{eq-nablc1}
\begin{array}{c}
\displaystyle
\nabla_{c'(0)}c'(0) = c''(0)+c'(0)(\Theta(c(0),c'(0))-\Theta(c(0),-c'(0)))=\\[10pt]
\displaystyle =\left(\frac{2}{3}f'''(0),-\frac{1}{2}\right)+O(e^{-2r}) \\[10pt]
\end{array}
\end{equation}
Using \eqref{eq-x2}, \eqref{eq-nablc1} we get
$$
\frac{\mathbf{g}_{c'(0)}(\nabla_{c'(0)}c'(0),c'(0))}{\mathrm{F}(c(0),c'(0))^2}=\frac{\frac{2}{3}f'''(0)g_{11}-\frac{1}{2}g_{12}}{\mathrm{F}(0,C e^{-2r}+O(e^{-3r}),1+O(e^{-2r}),0)^2} $$
Here $g_{ij}$ are calculated at the point $(0,C e^{-2r}+O(e^{-3r}),1+O(e^{-2r}),0)$. After substituting the values from \eqref{eq-F}, \eqref{eq-g11} we obtain
$$\frac{\mathbf{g}_{c'(0)}(\nabla_{c'(0)}c'(0),c'(0))}{\mathrm{F}(c(0),c'(0))^2}=-\frac{f'''(0)}{3}+O(e^{-2r})$$
Therefore,
\begin{equation}\label{eq-nablacdot}\nabla_{\dot{c}(0)}\dot{c}(0)=\frac{\left(f'''(0),-\frac{1}{2}\right)+(1,1)O(e^{-2r})}{\mathrm{F}(c(0),c'(0))^2}\end{equation}
Taking into account \eqref{eq-F},
$$\nabla_{\dot{c}(0)}\dot{c}(0)=\left(4f'''(0),-2\right)e^{-2r}+(1,1)O(e^{-3r})$$

Calculate the Rund curvature \eqref{eq-curvr} using the formulae \eqref{eq-x2}, \eqref{eq-nablacdot}.

$$
\begin{array}{c}
\displaystyle
\mathbf{k}_R(r)^2=\mathrm{F}(c(0),\nabla_{\dot{c}(0)}\dot{c}(0))=\\[10pt]
\displaystyle =\frac{\mathrm{F}\left(0,C e^{-2r}+O(e^{-3r}),-f'''(0)+O(e^{-2r}),\frac{1}{2}+O(e^{-2r})\right)}{\mathrm{F}(0,C e^{-2r}+O(e^{-3r}),1-Ce^{-2r}+O(e^{-3r}),0)^2} \\[10pt]
\end{array}
$$

From \eqref{eq-F}, \eqref{eq-Fl} we get
\begin{equation}\label{eq-formrund}\mathbf{k}_R(r)^2=1+C\left(\frac{2}{L}-\frac{8f'''(0)^2}{9}\right)e^{-2r}+O(e^{-3r})\end{equation}
Here $L>0$ is the length of the chord $\ell$ of $\partial U$ in the direction $(f'''(0),-1/2)$. Proposition 1 gives the uniform bounds on the curvature of $\partial U$ and proposition 2 claims that the angle between the chord $\ell$ and $x_2$ is uniformly separated from $\pi/2$, thus we conclude that $\frac{2}{L}$ is bounded from above.

Calculate the Finsler curvature  \eqref{eq-curvf} using the formulae \eqref{eq-x2}, \eqref{eq-nablacdot}.

$$
\begin{array}{c}
\displaystyle
\mathbf{k}_F(r)^2=\mathbf{g}_{\dot{c}(0)}(\nabla_{\dot{c}(0)}\dot{c}(0),\nabla_{\dot{c}(0)}\dot{c}(0))=\\[10pt]
\displaystyle =\frac{f'''(0)^2g_{11}-f'''(0)g_{12}+\frac{1}{4}g_{22}}{\mathrm{F}(0,C e^{-2r}+O(e^{-3r}),1-Ce^{-2r}+O(e^{-3r}),0)^4}, \\[10pt]
\end{array}
$$
Here $g_{ij}$ are considered at the point $(0,C e^{-2r}+O(e^{-3r}),1+O(e^{-2r}),0)$. Finally, from \eqref{eq-F},  \eqref{eq-g11} we obtain that
\begin{equation}\label{eq-formfins}\mathbf{k}_F(r)^2=1+C\left(-\frac{8}{9}f'''(0)^2+4f^{(4)}(0)\right)e^{-2r}+O(e^{-3r})\end{equation}

Proposition 2 gives the uniform bounds on the derivatives of $f$, and Theorem 1 is proved.

Note that the normal curvature $\mathbf{g}_\mathbf{n}(\nabla_{\dot{c}(s)}\dot{c}(s),\mathbf{n})$ of a hypersurface at the point $x$  depends only on the tangent vector to the curve $c(s)$ at $x$ (\cite{Sh}).
So, in order to obtain the normal curvature of the Hilbert hypersphere $S_r$ centred at $o$ at the point $p$ in the tangent direction $w$, we consider the normal curvature of the circle $S_r \cap \Pi$ which lies in the plane $\Pi = span (w,\vec{op})$.

From \eqref{eq-nablc1} we get the normal curvature of the circle of radius $r$.
\begin{equation}\label{eq-formnorm1}\mathbf{k}_{\mathbf{n}}(r)=\mathbf{g}_{\mathbf{n}}(\nabla_{\dot{c}(0)}\dot{c}(0),\mathbf{n})=\frac{\mathbf{g}_{\mathbf{n}}(c''(0),\mathbf{n})}{\mathrm{F}(c(0),c'(0))^2}\end{equation}
From the equality $g_{12}(0,x_2,0,1)=0$ (\eqref{eq-garb}) it follows that the unit normal vector $\mathbf{n}$ to the circle at $(0,x_2)$ is exactly $\frac{1}{\mathrm{F}(0,x_2,0,1)}(0,-1)$.

Finally, getting into account \eqref{eq-F}, \eqref{eq-x2}, \eqref{eq-c11}, \eqref{eq-F2}, \eqref{eq-garb}:
\begin{equation}
\label{eq-formnorm}
\begin{array}{c}
\displaystyle \mathbf{k}_{\mathbf{n}}(r)=\frac{\frac{1}{2}g_{22}(0,C e^{-2r}+O(e^{-3r}),0,1)}{\mathrm{F}(0,C e^{-2r}+O(e^{-3r}),1-Ce^{-2r}+O(e^{-3r}),0)^2 \mathrm{F}(0,C e^{-2r}+O(e^{-3r}),0,1)}=\\[10pt]
\displaystyle  =1+C\left(\frac{1}{H}-\frac{8f'''(0)^2}{9}\right)e^{-2r}+O(e^{-3r})\\[10pt]
\end{array}
\end{equation}

If the Euclidean normal curvatures of the hypersurface $\partial U$ are bounded ($k_2\leqslant k_n \leqslant k_1$) then the curvature of the curve  $\partial U' = \partial U \cap \Pi$ is bounded as well. Indeed consider the point $x \in \partial U' \subset \partial U$. Then the curvature $k(x)$ of $\partial U'$ and the normal curvature $k_n(x)$ of $\partial U$ are related as $k(x) = \frac{k_n(x)}{\cos \beta}$. Here $\beta$ in the angle between the radial and normal direction to $\partial U$ at $x$. Using lemma \eqref{lem-angle} we find that $\omega_0 k_2\leqslant \cos \beta \leqslant 1$. Hence the curvature of $\partial U'$ is uniformly bounded for all $y$. Applying proposition 1 for the Hilbert geometry based on $U'$ we get the uniformity of the series expansion \eqref{eq-formnorm} which ends the proof of Theorem 2.


\begin{thebibliography}{99}

\bibitem{BCS}
Bao D., Chern S. S., Shen Z. An Introduction to Riemann-Finsler Geometry. -- Springer-Verlag, 2000.
\bibitem{Boo2}
Borisenko A.A. Convex sets in Hadamard manifolds, Differential Geometry and its Applications 17 (2002), 111-121.
\bibitem{Boo3}
Borisenko A. A., Olin E. A. Asymptotic Properties of Hilbert Geometry // Journal of Math. Phys., Anal., Geometry. -- 2008. --  Vol. 4, No 3.  -- P. 327-345.

\bibitem{CV}
Colbois B., Verovic. P. Rigidity of Hilbert Metrics // Bull. Austral. Math. Soc. -- 2002. --  No 65. -- P. 23-34.

\bibitem{F}
Finsler P. \"Uber Kurven und Fl\"achen in allgemeinen R\"aumen. -- Birkh\"auser Verlag -- Basel, 1951.

\bibitem{Sh}
Shen Z. Lectures on Finsler Geometry. -- World Scientific Publishing Co, 2001. -- 306 p.

\bibitem{RR}
Rund H. The Differential Geometry Of Finsler Spaces, 1959




\end{thebibliography}
 \end{document}